\documentclass[a4paper,10pt]{article}
\usepackage[english]{babel}
\usepackage{graphicx}
\usepackage[english]{babel}
\usepackage{amsmath}
\usepackage{amssymb}
\usepackage{amsfonts}
\usepackage{enumerate}
\usepackage{pstricks}
\usepackage{cite}
\usepackage{hyperref}
\usepackage{float}
\usepackage{comment}
\usepackage{authblk}

\usepackage{accents}

\usepackage[top=2cm, bottom=2cm, left=2cm, right=2cm]{geometry}

\newcommand*{\BAR}[1]{\overline{#1}}

\newtheorem{theorem}{Theorem}[section]
\newtheorem{thm}[theorem]{Theorem}
\newtheorem{prop}[theorem]{Proposition}
\newtheorem{cor}[theorem]{Corollary}
\newtheorem{defn}[theorem]{Definition}

\newtheorem{rem}[theorem]{Remark}

\newenvironment{proof}[1][Proof]{\textbf{#1: } }{\begin{flushright}$\blacksquare$\end{flushright}}

\begin{document}

\title{Three-Dimensional Analogues of the Blasius–Chaplygin Formulas}
\author{Dmitrii Legatiuk\thanks{Chair of Mathematics, Universit\"at Erfurt, Germany}, Heikki Orelma\thanks{Tampere Institute of Mathematics, Hepolamminkatu 51, 33720 Tampere, Finland
E-Mail: heikki.orelma@proton.me
}}
\date{}

\maketitle

\begin{abstract}
The classical Blasius--Chaplygin formula provides an elegant method for calculating the lift force on a two-dimensional body in steady, irrotational flow. The key ingredient is the definition of a complex-valued potential function 
\[
f(z) = \varphi(x, y) + \mathbf{i}\psi(x, y),
\]
which can then be integrated using Cauchy's theorem along any closed contour surrounding the body.

In this paper, we propose a three-dimensional extension of the classical Blasius--Chaplygin formula using quaternionic analysis. After presenting the basics of quaternionic analysis, we discuss how \emph{monogenic functions}---the quaternionic analog of classical holomorphic functions---can be used to describe problems in fluid dynamics. Finally, we present the Blasius--Chaplygin formula in quaternionic form.

\end{abstract}

\section{Introduction}

The velocity potential $\phi$ and stream function $\psi$ are used to describe the components of the velocity field in the case of the classical two-dimensional ideal (irrotational and inviscid) gas dynamics through the following expressions \cite{LL,LS}:
\begin{equation}\label{nop}
\begin{array}{ccl}
v_1 & = & \displaystyle \frac{\partial \phi}{\partial x}=\frac{\partial \psi}{\partial y},\\
v_2 & = & \displaystyle \frac{\partial \phi}{\partial y}=-\frac{\partial \psi}{\partial x}.
\end{array}
\end{equation}
Stream functions are geometrically related to the flow since streamlines are their level curves. Moreover, equations~(\ref{nop}) form the {\itshape Cauchy-Riemann system} known from complex analysis, which consequently leads to defining the complex potential function $f(\mathsf{z}) := \phi + \mathbf{i} \psi$. This complex function describes the flow’s potential and stream function fields such that its derivative $f'(\mathsf{z})$ corresponds to the local flow velocity. Through this, complex analysis opens an exceptionally efficient and elegant mathematical tool for fluid mechanics, especially for the study of ideal two-dimensional flows, see again the classical references \cite{LL,LS}.\par
A central problem in the ideal flow mechanics is a calculation of forces, which requires computing pressure forces arising from the distribution of flow velocity on the body’s surface. In the two-dimensional case, the complex potential provides a method for this, where the force acting on the body is given by
\begin{equation*}
F_x - \mathbf{i} F_y = \frac{\mathbf{i}\rho}{2} \oint_C \left( f'(\mathsf{z}) \right)^2\, d\mathsf{z}.
\end{equation*}
This formula is called the \textit{Blasius-Chaplygin formula}. It is noted that from a mathematical point of view, holomorphy of the potential inside the region bounded by the path $C$ implies that the flow exerts no forces on the body. Consequently, forces appear only when the holomorphy is lost at some point of the body. Complex analysis provides many tools to analyze such situations.\par
The complex integral form of the force formula known as the Blasius–Chaplygin formula was developed independently by two researchers in the early 20th century. Let us now provide a brief overview of these discoveries:
\begin{itemize}
\item {\itshape Paul Blasius} (1883–1970), a German physicist and student of Prandtl, published the formula in 1911 in his article "Grenzschichten in Fl\"ussigkeiten mit kleiner Reibung" (Boundary layers in fluids with small viscosity) \cite{blasius1911}. Although the article is mostly known for Blasius’ boundary layer solution, he also used the complex flow theory for force calculations in a two-dimensional potential flow, presenting an integral formula to determine lift and drag on a body.
\item {\itshape Sergei A. Chaplygin} (1869–1942), a Russian mathematician and mechanics researcher, published similar formulas almost simultaneously (in 1902 and more extensively in 1904) as part of his studies on flows over lifting surfaces. His work appeared in the publications of the Russian Academy of Sciences under titles such as "On gas streams with uniform motion" \cite{chaplygin1902}. Chaplygin applied complex functions to analyze the flow and calculated forces on bodies using a similar integral method.
\end{itemize}\par
Although P. Blasius and S.A. Chaplygin mainly studied forces, a formula for the moment (about the origin) has also been derived based on their works:
\begin{equation*}
M = -\frac{\rho}{2} \operatorname{Re} \oint_C \mathsf{z} \left( f'(\mathsf{z}) \right)^2\, d\mathsf{z}.
\end{equation*}
Thus, the force and moment formulas combine complex analysis and potential flow theory, forming an elegant method to calculate the force and moment acting on a body in a two-dimensional inviscid flow.\par
It is important to emphasize that the Blasius formulas are derived for a two-dimensional ideal flow where the body considered is essentially infinitely long in the third dimension. Consequently, the flow is homogeneous in the third dimension and unaffected along the length. Therefore, the formulas give the values of flow force and moment per cross-sectional length, meaning they are {\itshape force and moment densities}:
\begin{itemize}
\item The unit of force is therefore Newtons per meter $[N/m]$, which means the force acting on a one-meter long portion of the body.
\item The moment is similarly expressed as moment per unit length, i.e., in units of Newton-meters per meter $[Nm/m]$, which typically simplifies to Newton-meters $[Nm]$ in the two-dimensional case.
\end{itemize}
In practice, this means that the Blasius formulas apply to describing, for example, the lift and torque generated by the cross-section of an airplane wing, assuming the wing is long and the flow is homogeneous in two dimensions.\par
This article proposes a three-dimensional counterpart of the complex potential and related Blasius-Chaplygin formulas. This three-dimensional extension has become possible with the advanced in the field of {\itshape quaternion analysis}, which is a hypercomplex complex function theory in $\mathbb{R}^{3}$ and $\mathbb{R}^{4}$, see \cite{Guerlebeck_1,Guerlebeck_2} and references therein for an overview of a hypercomplex analysis. The core idea of quaternionic analysis is to work with the so-called {\itshape monogenic functions}, which generalize the concept of holomorphic functions to quaternions. In the context of boundary value problems in mechanics, two general strategies of using tools of quaternionic analysis can be distinguished: first strategy is to apply the so-called {\itshape quaternionic operator calculus}, which allows to construct explicit representation formulas for the solution of boundary value problems by using singular integral operators, such as Cauchy-Bitsadze operator and Teodorescu transform, see \cite{Guerlebeck_3,Guerlebeck_4} for examples of using this strategy for elasticity theories; the second strategy is to generalize the classical analytical solution approaches to the quaternionic case, see for example works \cite{Bock_1,Bock_2} for spatial analogues of the classical Kolosov-Muskhelishvili formulas. In this paper, we follow the second strategy.\par
While two-dimensional potential flows are well understood through complex analysis, the extension to three dimensions is far from trivial. Quaternion analysis offers promising tools for this extension, enabling novel analytical methods for ideal flow problems in 3D. First ideas towards developing quaternionic theory for three-dimensional flow and the construction of monogenic flow potentials have been presented in \cite{O1,O2}. This paper presents a continuation of these works. Moreover, to make the paper as self-contained and complete as possible, we recall also some results from \cite{O1,O2} to support the reader.\par

\section{Preliminaries}

In this section, we provide a brief introduction into quaternions and quaternionic analysis. We refer to \cite{Guerlebeck_1,Guerlebeck_2} and references therein for more detailed information, as well as for applications of quaternionic setting.

\subsection{Quaternions}

Quaternions are a generalization of complex numbers $\mathbb{C}$ to $\mathbb{R}^{4}$. The algebra of quaternions is defined by the basis elements $\{ 1,\mathbf{i},\mathbf{j},\mathbf{k}\}$, which satisfy the well-known multiplication rules
\begin{equation*}
\mathbf{i}^2=\mathbf{j}^2=\mathbf{k}^2=\mathbf{i}\mathbf{j}\mathbf{k}=-1.
\end{equation*}
The algebra of quaternions is associative, but not commutative, as illustrated by the following multiplication rules for the unit elements
\begin{equation*}
\begin{array}{rclrclrcl}
\mathbf{i}\mathbf{j} & = & \mathbf{k}, & \mathbf{j}\mathbf{k} & = & \mathbf{i}, & \mathbf{k}\mathbf{i} & = & \mathbf{j}, \\
\mathbf{j}\mathbf{i} & = & -\mathbf{k}, & \mathbf{k}\mathbf{j} & = & -\mathbf{i}, & \mathbf{i}\mathbf{k} & = & -\mathbf{j}.
\end{array}
\end{equation*}
The element $1$ is the multiplicative identity in the algebra of quaternions.\par
A general quaternion $q\in\mathbb{H}$ is of the form
\begin{equation*}
q=q_0+q_1\mathbf{i}+q_2\mathbf{j}+q_3\mathbf{k},
\end{equation*}
where $q_0,q_1,q_2,q_3\in\mathbb{R}$. Since a complex number $z\in\mathbb{C}$ is of the form $z=x+iy$, a quaternion can be viewed as a “complex number” with three imaginary units instead of one. Geometrically, quaternions can be identified with a four-dimensional vector space, that is, $\mathbb{H}\cong\mathbb{R}^4$, where the vector space isomorphism is defined as
\begin{equation*}
q_0+q_1\mathbf{i}+q_2\mathbf{j}+q_3\mathbf{k} \quad \mapsto \quad (q_0,q_1,q_2,q_3).
\end{equation*}
This identification holds only when quaternions are regarded as a vector space, since the space $\mathbb{R}^4$ does not have a naturally defined multiplication.\par
For convenience reasons, we will sometimes switch to the following denotation of the basis elements:
\begin{equation}\label{emuuttujat}
\mathbf{e}_0 :=1, \quad \mathbf{e}_1 :=\mathbf{i}, \quad \mathbf{e}_2 := \mathbf{j}, \quad \mathbf{e}_3 := \mathbf{k}.    
\end{equation}\par
Further, a quaternion can be written in the form
\begin{equation*}
q=q_0+\mathbf{q},
\end{equation*}
where $q_0$ is the \textit{scalar part} of the quaternion $q$ and $\mathbf{q}=q_{1}\mathbf{e}_{1}+q_{2}\mathbf{e}_{2}+q_{3}\mathbf{e}_{3}$ is the \textit{vector part}, which are denoted as $\text{Sc}(q):=q_0$ and $\text{Vec}(q):=\mathbf{q}$, respectively. The \textit{conjugate} of a quaternion is defined by
\begin{equation*}
\overline{q} := q_0-\mathbf{q}. 
\end{equation*}
Quaternion conjugation reverses the order of multiplication, that is,
\begin{equation*}
\overline{pq}=\overline{q}\,\overline{p}, \qquad p,q\in\mathbb{H}.
\end{equation*}
The scalar and vector parts of a quaternion $q\in\mathbb{H}$ can be computed using the conjugate as
\[
\text{Sc}(q)=\frac{1}{2}(q+\overline{q}),\ \text{ and }\ \text{Vec}(q)=\frac{1}{2}(q-\overline{q}).
\]
Additionally, we have the following identity
\begin{equation*}
\overline{q}q=q\overline{q}=q_0^2+q_1^2+q_2^2+q_3^2.
\end{equation*}
This means that the Euclidean norm can be defined by setting
\begin{equation*}
|q|^2=\overline{q}q=q\overline{q}.
\end{equation*}
Consequently, the multiplicativity of the norm with respect to the product follows immediately:
\begin{equation*}
|pq|=|p||q|.
\end{equation*}
The Euclidean inner product induces the Euclidean norm, and it can be expressed using the quaternionic product as
\begin{equation*}
\langle q,p \rangle =\text{Sc}(q\overline{p})=\sum_{n=0}^3 p_{n}q_{n},
\end{equation*}
for $p,q\in\mathbb{H}$. Moreover, the scalar part of the product satisfies the symmetry relations
\begin{equation*}
\text{Sc}(q\overline{p})=\text{Sc}(p\overline{q})=\text{Sc}(\overline{q}p)=\text{Sc}(\overline{p}q).
\end{equation*}

\subsection{Quaternionic analysis}

Although there are several alternative approaches to develop a quaternionic analysis, see again \cite{Guerlebeck_1,Guerlebeck_2}, the underlying idea is always to construct a class of functions analogous to analytic functions in complex analysis. Keeping in mind applications to fluid mechanics, let us consider functions of the form
\begin{equation*}
f\colon U\to\mathbb{H},
\end{equation*}
where $U\subset \mathbb{R}^3$. Accordingly, the functions are of the form
\begin{equation*}
f(x,y,z)=f_0(x,y,z)+f_1(x,y,z)\mathbf{i}+f_2(x,y,z)\mathbf{j}+f_3(x,y,z)\mathbf{k},
\end{equation*}
where the component functions $f_0,f_1,f_2,f_3\colon U\to\mathbb{R}$. The usual properties of quaternion-valued functions -- such as differentiability, boundedness, integrability etc. -- are defined in terms of the corresponding properties of the component functions. In this article, the points $(x,y,z)$ are identified with the \textit{reduced quaternions} 
\begin{equation*}
\mathsf{x}=x+y\mathbf{i}+z\mathbf{j}.
\end{equation*}\par
Next, we define a class of functions in quaternionic analysis, which generalize the concept of complex analytic functions to $\mathbb{R}^{3}$. First, we need to introduce the the \textit{generalized Cauchy-Riemann operator}:
\begin{equation*}
D := \frac{\partial}{\partial x}+\mathbf{i}\frac{\partial}{\partial y}+\mathbf{j}\frac{\partial}{\partial z}
\end{equation*}
and its conjugate
\begin{equation*}
\overline{D}:=\frac{\partial}{\partial x}-\mathbf{i}\frac{\partial}{\partial y}-\mathbf{j}\frac{\partial}{\partial z}.
\end{equation*}
These operators give rise to the definition:
\begin{defn}[Monogenicity in a domain $U$]
Let $U\subset\mathbb{R}^3$ be an open set and $f\colon U\to\mathbb{H}$ a differentiable function. The function is called left-monogenic if $Df=0$ in $U$, and correspondingly right-monogenic if $fD=0$ in $U$. 
\end{defn}\par
As it has been already mentioned, the monogenicity concept introduced in this definition corresponds to the complex analyticity. The following proposition explores this connection further by introduce system of partial differental equations analogous to the Cauchy-Riemann equations:
\begin{prop}
Let $f\colon U\to\mathbb{H}$ be a differentiable function. The following conditions are equivalent:
\begin{itemize}
\item[a)] $Df=0$ ($f$ is monogenic),
\item[b)] The components of the function $f=f_0+f_1\mathbf{i}+f_2\mathbf{j}+f_3\mathbf{k}$ are a solution of the \textit{Moisil-Theodorescu system} 
\begin{equation*}
\left\{\begin{array}{rcl}
\displaystyle \frac{\partial f_0}{\partial x}-\frac{\partial f_1}{\partial y}-\frac{\partial f_2}{\partial z} & = & 0,\\
\\
\displaystyle \frac{\partial f_1}{\partial x}+\frac{\partial f_0}{\partial y}+\frac{\partial f_3}{\partial z} & = & 0,\\
\\
\displaystyle \frac{\partial f_2}{\partial x}+\frac{\partial f_0}{\partial z}-\frac{\partial f_3}{\partial y} & = & 0,\\
\\
\displaystyle \frac{\partial f_3}{\partial x}+\frac{\partial f_2}{\partial y}-\frac{\partial f_1}{\partial z} & = & 0.
\end{array} \right.
\end{equation*}
\end{itemize}
\end{prop}\par
The principal difference between monogenic functions in quaternionic analysis and analytic functions in complex analysis is that, unlike the powers of complex numbers $(x+iy)^n$, the powers of the quaternionic variable $(x+\mathbf{i}y+\mathbf{j}z)^n$ are no longer monogenic. This can be seen, for example, by computing $D(x+\mathbf{i}y+\mathbf{j}z)=-1$. In quaternionic analysis, much attention has been devoted to constructing various bases of monogenic polynomials, see for example works and references therein \cite{Bock_3,Cruz,Grigoriev}.\par
The definition of monogenicity can also be associated with the conjugate Cauchy–Riemann operator. In this case, instead of monogenicity, one speaks of antimonogenicity.
\begin{defn}[Antimonogenicity in a domain $U$]
Let $U\subset\mathbb{R}^3$ be an open set and $f\colon U\to\mathbb{H}$ a differentiable function. The function is called left-antimonogenic if $\BAR{D}f=0$ in $U$, and correspondingly right-antimonogenic if $f\BAR{D}=0$ in $U$. 
\end{defn}\par
One of the most important property of the generalized Cauchy–Riemann operator and its conjugate is the factorization of the \textit{Laplace operator}:
\begin{equation*}
D\overline{D} = \overline{D}D = \Delta \mbox{ with } \Delta := \frac{\partial^2}{\partial x^2}+\frac{\partial^2}{\partial y^2}+\frac{\partial^2}{\partial z^2}.
\end{equation*}
Consequently, the component functions of every monogenic function are harmonic. Conversely, monogenic functions can be generated from harmonic functions. Let us now consider two methods to generate a monogenic function from a harmonic function. The first method follows directly from the factorization property:
\begin{prop}
Let $u\colon \Omega\to\mathbb{R}$ be a harmonic function in an open set $\Omega\subset \mathbb{R}^3$. Then $\overline{D}u$ is a monogenic function in the set $\Omega$. 
\end{prop}\par
The second method, which is more useful in fluid dynamics, is based on the formula presented in \cite{BDS}, where the case $n=2$ corresponds to quaternionic analysis. To make the paper self-containing, we also provide a proof of this result, since the original source does not present a proof of the formula either.
 The proof uses the \textit{Euler operator} defined as follows
\begin{equation*}
E := x\frac{\partial }{\partial x}+y\frac{\partial }{\partial y}+z\frac{\partial }{\partial z},
\end{equation*}
which can be expressed using the generalized Cauchy–Riemann operator and inner product as $E=\langle \mathsf{x},D\rangle$. Hence, we have the following proposition:
\begin{prop}\label{HarmRealOs}
Let $\Omega\subset\mathbb{R}^3$ be an open set star-shaped with respect to the origin and let $u\colon\Omega\to\mathbb{R}$ be a harmonic function. Then
\begin{equation*}
f(\mathsf{x})=u(\mathsf{x})+\text{Vec}\left(\int_0^1  \overline{D}u(\mathsf{x}t)\mathsf{x}t\; dt\right)
\end{equation*}
is a monogenic function whose real part is the function $u$. 
\end{prop}
\begin{proof}
To make construction conciser, we use the notation~(\ref{emuuttujat}) for the unit elements along with the notations $x_0 := x$, $x_1 := y$, and $x_2 := z$.\par 
We start by writing the function $f$ as follows
\begin{equation*}
f(\mathsf{x})=u(\mathsf{x})+\int_0^1  \overline{D}u(\mathsf{x}t)\mathsf{x}t\; dt-\text{Sc}\left(\int_0^1  \overline{D}u(\mathsf{x}t)\mathsf{x}t\; dt\right).
\end{equation*}
Let us now introduce the notation
\begin{equation*}
I(\mathsf{x}) := \int_0^1  \overline{D}u(\mathsf{x}t)\mathsf{x}t\; dt,
\end{equation*}
implying that the function $f$ can now be rewritten as
\begin{equation*}
f(\mathsf{x})=u(\mathsf{x})+I(\mathsf{x})-\text{Sc}(I(\mathsf{x})).
\end{equation*}
By applying the fact that $D\overline{D}u=\Delta u=0$ and the characterization of reduced quaternions\footnote{A quaternion $p\in\mathbb{H}$ is a reduced quaternion (i.e., of the form $q=q_0\mathbf{e}_{0}+q_1\mathbf{e}_1+q_2\mathbf{e}_2$) if and only if
\begin{equation*}
\sum_{j=0}^{2}\mathbf{e}_{j}q\mathbf{e}_{j} = -\overline{q}.
\end{equation*}
}, we obtain
\begin{equation*}
\begin{array}{rcl}
DI(\mathsf{x}) & = & \displaystyle \int_0^1  D\overline{D}u(\mathsf{x}t)\mathsf{x}t\; dt+\int_0^1  \sum_{j=0}^2 \mathbf{e}_j\overline{D} u(\mathsf{x}t)\frac{\partial \mathsf{x}}{\partial x_j}t\; dt\\
& = & \displaystyle \int_0^1 \sum_{j=0}^2 \mathbf{e}_j\overline{D}u(\mathsf{x}t)\mathbf{e}_{j}t\; dt =- \int_0^1  Du(\mathsf{x}t)t\; dt.
\end{array}
\end{equation*}
Additionally, we have the following relation
\begin{equation*}
\text{Sc}(\overline{D}u(\mathsf{x}t)\mathsf{x})=\langle Du(\mathsf{x}t),\mathsf{x}\rangle=\sum_{j=0}^2 tx_j\frac{\partial }{\partial (t x_j) }u(\mathsf{x}t)=Eu(\mathsf{x}t).
\end{equation*}
Thus,
\begin{equation*}
Df(\mathsf{x})=Du(\mathsf{x})- \int_0^1  Du(\mathsf{x}t)t\; dt-\int_0^1  DEu(\mathsf{x}t)t\; dt.
\end{equation*}
A direct computation yields
\begin{equation*}
DEu=Du+EDu.
\end{equation*}
Hence,
\begin{equation*}
Df(\mathsf{x}) = Du(\mathsf{x})-2 \int_0^1  Du(\mathsf{x}t)t\; dt-\int_0^1  EDu(\mathsf{x}t)t\; dt.
\end{equation*}
Now, $EDu(t\mathsf{x})=t\dfrac{d}{dt}Du(\mathsf{x}t)$, so by integration by parts we get
\begin{equation*}
\int_0^1  EDu(\mathsf{x}t)t\; dt=\int_0^1  t^2\frac{d}{dt}Du(\mathsf{x}t)\; dt=[t^2Du(\mathsf{x}t)]^1_0-2\int_0^1  t Du(\mathsf{x}t)\; dt 
\end{equation*}
where $[t^2Du(\mathsf{x}t)]^1_0=Du(\mathsf{x})$.
Therefore,
\begin{equation*}
Df(\mathsf{x})=Du(\mathsf{x})-2 \int_0^1  Du(\mathsf{x}t)t\; dt-Du(\mathsf{x})+2\int_0^1  tDu(\mathsf{x}t)\; dt=0. 
\end{equation*}
\end{proof}\par
The next step is to recall basics of integration theory for monogenic functions, where the integration has to be computed now over surfaces, and not just curves as in the two-dimensional case. Assume that $\Omega\subset\mathbb{R}^3$ is a compact set. Suppose that the boundary surface $\partial \Omega$ of the set $\Omega$ is piecewise smooth, meaning that an outward unit normal can be defined almost everywhere on the surface. Such a set $\Omega$ will henceforth be called \textit{regular}. We define the quaternion-valued $2$-differential form
\begin{equation}\label{pintaalkio}
 d\sigma =dy\wedge dz-\mathbf{i}dx\wedge dz+\mathbf{j}dx\wedge dy,   
\end{equation}
and assume that $n=n_1+\mathbf{i}n_2+\mathbf{j}n_3$ is a quaternionic representation of the outward unit normal on the surface $\partial \Omega$. It can be shown\footnote{The proof is straightforward if one assumes that the surface has a smooth parametrization on some chart $r(u,v)=(x(u,v),y(u,v),z(u,v))$. One then simply restricts $d\sigma$ to this chart and observes that $n(r(u,v))=\pm \frac{r_u(u,v)\times r_v(u,v)}{|r_u(u,v)\times r_v(u,v)|}$. In this case, $dS=\pm |r_u(u,v)\times r_v(u,v)|\;  du\wedge dv$. The plus-minus combinations above are chosen so that the unit normal vector points outward.} that
\begin{equation}\label{normaali}
d\sigma= n\; dS
\end{equation}
when restricted to the surface $\partial \Omega$, where $dS$ is the scalar surface element.\par
The significance of the above quaternion-valued surface element $d\sigma$ in quaternionic analysis is revealed in the following integral formula. The proof of the result is a "one-line proof" based on the classical Stokes' theorem.
\begin{thm}[Quaternionic Stokes' theorem]\label{STOKES}
Let $\Omega$ be a regular domain and let $f,g\colon \Omega\to \mathbb{H}$ be continuous functions such that $f,g\colon\Omega^\circ\to \mathbb{H}$ are differentiable, where $\Omega^\circ$ denotes the interior points of $\Omega$. Then
\begin{equation*}
\int_{\partial \Omega} g d\sigma f=\int_\Omega\big((gD)f+g(Df)\big)dV,
\end{equation*}
where $dV=dx\wedge dy\wedge dz$.
\end{thm}\par
An important example of a function that is monogenic on both the left and the right in the integration theory is the \textit{Cauchy kernel}, defined as follows
\begin{equation*}
\mathcal{E}(\mathsf{x})=\frac{1}{4\pi} \frac{\overline{\mathsf{x}}}{|\mathsf{x}|^3},
\end{equation*}
where $\mathsf{x}=x+y\mathbf{i}+z\mathbf{j}$. A direct computation shows that $D\mathcal{E}(\mathsf{x})=\mathcal{E}(\mathsf{x})D=0$ whenever $\mathsf{x}\neq 0$. By substituting the fundamental solution for $g$ in the Stokes' formula, and assuming monogenicity of $f$, we immediately obtain the Cauchy formula:
\begin{thm}[Cauchy Integral Formula]
Let $\Omega$ be a regular domain and let $f\colon\Omega\to \mathbb{H}$ be continuous such that $f\colon\Omega^\circ\to \mathbb{H}$ is monogenic. Then
\begin{equation*}
f(\mathsf{x})=\int_{\partial\Omega} \mathcal{E}(\mathsf{y}-\mathsf{x}) d\sigma(\mathsf{y}) f(\mathsf{y}).
\end{equation*}
\end{thm}\par
In flow dynamics, however, we largely work with surface integrals. For this reason, the quaternionic Stokes' theorem is particularly important, as it connects monogenicity with surface integrals. Hence, we obtain the following result:
\begin{prop}\label{ZIL}
Let $\Omega$ be a regular domain and let $f,g\colon\Omega\to \mathbb{H}$ be continuous such that $g\colon\Omega^\circ\to \mathbb{H}$ is right-monogenic and $f\colon\Omega^\circ\to \mathbb{H}$ is left-monogenic. Then
\begin{equation*}
\int_{\partial \Omega} g d\sigma f=0.
\end{equation*}
\end{prop}
In order to make use of the monogenicity concept, for example in the form of a three-dimensional residue calculus, all surface integrals must be expressed, in one way or another, in terms of integrals of the type shown above.\par

\section{Monogenic flow potential}

In this section, we present the use of monogenic functions to introduce a flow potential. Let us consider a velocity vector field of an ideal flow $\vec{v} = (v_1, v_2, v_3)$ defined in an open set $U \subset \mathbb{R}^3$, and assume that $\phi \colon U \to \mathbb{R}$ is a \textit{potential of the velocity vector field}, that is,
\begin{equation*}
v_1 = \frac{\partial \phi}{\partial x}, \quad v_2 = \frac{\partial \phi}{\partial y}, \quad v_3 = \frac{\partial \phi}{\partial z}.
\end{equation*}
Further, we can consider the \textit{stream function} is $\Psi = \psi_1\mathbf{i}+\psi_2\mathbf{j}+\psi_3\mathbf{k}$ such that 
\begin{equation*}
w=\phi+\Psi
\end{equation*}
is a monogenic function in $U$. The function $w$ is called a \textit{monogenic flow potential}. Because
\begin{equation*}
v:=v_1+v_2\mathbf{i}+v_3\mathbf{j} = D\phi,
\end{equation*}
and $Dw=0$, then
\begin{equation*}
D\Psi=-v_1-v_2\mathbf{i}-v_3\mathbf{j}. 
\end{equation*}
From this it follows that
\begin{equation} \label{Nopeusv}
 v=v_1+v_2\mathbf{i}+v_3\mathbf{j}=\frac{1}{2}D\overline{w}.   
\end{equation}
Note that since the components of the monogenic flow potential are harmonic functions, in ideal flow, the velocity is a left anti-monogenic function
\[
\overline{D}v=0.
\]
A monogenic flow potential is a function, which is constructed starting with a harmonic velocity potential $\phi$.\par 
Proposition~\ref{HarmRealOs} leads to the following result:
\begin{prop}
Let $U\subset\mathbb{R}^3$ be a star-shaped open set with respect to the origin. Then a monogenic flow potential exists in this domain. 
\end{prop}
In domains more general than star-shaped domains, the construction of a monogenic flow potential is an open question. Moreover, the stream function, and thus the monogenic flow potential, is not unique, because it is always possible to add a monogenic vector-valued function to the stream function $\Psi$. More precisely, the transformation $\Psi \mapsto \Psi + \mathbf{H}$ is a gauge transformation that leaves the flow properties unchanged. In particular, the vector-valued monogenic function $\mathbf{H}$ does not depend on the velocity potential $\phi$.\par
In general, it is difficult to geometrically interpret the stream function $\Psi$. Nonetheless, the following theorem presents a monogenic flow potential, which produces a stream functions, that can be easily interpreted geometrically:
\begin{thm}[\cite{O2}]\label{GeomFlowPT}
 If the monogenic flow potential is of the form   
\begin{equation*}
w(x,y,z)=\phi(x,y,z) + \psi_1(x,y)\mathbf{i} + \psi_2(x,z)\mathbf{j} + \psi_3(y,z)\mathbf{k},
\end{equation*}
satisfying the conditions 
\begin{equation*}
\begin{array}{rclrclrcl}
\displaystyle \frac{\partial \psi_1}{\partial y} & = & \displaystyle \frac{1}{2}v_1, & \displaystyle \frac{\partial \psi_2}{\partial z} & = & \displaystyle \frac{1}{2}v_1, & 
\displaystyle \frac{\partial \psi_3}{\partial z} & = & \displaystyle -\frac{1}{2} v_2, \\
\\
\displaystyle \frac{\partial \psi_1}{\partial x} & = & \displaystyle -\frac{1}{2}v_2, & \displaystyle \frac{\partial \psi_2}{\partial x} & = & \displaystyle -\frac{1}{2} v_3, &
\displaystyle \frac{\partial \psi_3}{\partial y} & = & \displaystyle \frac{1}{2} v_3,
\end{array}
\end{equation*} 
then
\begin{equation*}
\nabla\phi\cdot \nabla\psi_1=\nabla\phi\cdot \nabla\psi_2=\nabla\phi\cdot \nabla\psi_3=0
\end{equation*}
and
\begin{equation*}
v_1 = \frac{\partial \psi_1}{\partial y}+\frac{\partial \psi_2}{\partial z}, \quad v_2 = -\frac{\partial \psi_1}{\partial x}-\frac{\partial \psi_3}{\partial z}, \quad v_3 = -\frac{\partial \psi_2}{\partial x}+\frac{\partial \psi_3}{\partial y}.
\end{equation*}
\end{thm}
The idea of the above stream function is to search for two-dimensional stream functions on the coordinate planes. Accordingly, the streamlines lie on the intersection of the level surfaces defined by them. We refer to \cite{O2} for further details.\par

\section{Blasius–Chaplygin formula for the force}

Let $U$ be a body that is assumed to form a regular domain, and let $p(x,y,z)$ be the pressure field on its surface $\partial U$. Then the total force caused by the pressure
\begin{equation*}
F=F_1+F_2\mathbf{i}+F_3\mathbf{j}
\end{equation*}
can be calculated using the formula \cite{LL}:
\begin{equation*}
F=-\int_{\partial U} p\, d\sigma.
\end{equation*}\par
Let us recall the famous \textit{Bernoulli equation} along each streamline, neglecting the effect of gravity, see e.g. \cite{LL}), which is given by
\begin{equation*}
p+\frac{\rho}{2}|v|^2=\text{const.}
\end{equation*}
Note that in the case of an ideal flow, the density $\rho$ is constant. Substituting the Bernoulli equation into the expression for the force yields
\begin{equation*}
F=\frac{\rho}{2}\int_{\partial U} |v|^2 \, d\sigma.
\end{equation*}
Above, the integral over the constant term vanishes by Theorem~\ref{STOKES}. Substituting the expression for the potential from~(\ref{Nopeusv}), we obtain
\begin{equation}\label{BlasiusVoima}
F=\frac{\rho}{8}\int_{\partial U} |D\overline{w}|^2 \, d\sigma.
\end{equation}
The above expression is a spatial analog of the \textit{Blasius–Chaplygin force formula}, which allows the force acting on the body to be calculated using the monogenic flow potential. However, the previous integral is not of the form presented in the integration theory, since $|D\overline{w}|^2$ does not represent a monogenic function when $w$ is monogenic. Nevertheless, by taking the conjugate of the force expression on both sides and extracting the scalar part, we obtain
\begin{equation*}
F_1 = \frac{\rho}{8}\text{Sc}\int_{\partial U} |D\overline{w}|^2 \, \overline{d\sigma} = \frac{\rho}{8}\text{Sc}\int_{\partial U}  (w\overline{D}) (D\overline{w}) \overline{d\sigma}  
\end{equation*}
and by multiplying both sides by $\mathbf{i}$, we obtain
\begin{equation*}
F_2 = \frac{\rho}{8}\text{Sc}\int_{\partial U} |D\overline{w}|^2 \, d\sigma \mathbf{i} = \frac{\rho}{8}\text{Sc}\int_{\partial U}  (w\overline{D}) (D\overline{w}) \overline{d\sigma} \mathbf{i}.
\end{equation*}
Similarly to the above, we obtain
\begin{equation*}
F_3 = \frac{\rho}{8}\text{Sc}\int_{\partial U}  (w\overline{D}) (D\overline{w}) \overline{d\sigma} \mathbf{j}.
\end{equation*}
Assume that $w$ is of the same form as in Theorem~\ref{GeomFlowPT}. Define a quaternion-valued $1$-form by setting
\begin{equation*}
U=2(-\psi_3+\phi \mathbf{k})dx+2(-\psi_2-\phi \mathbf{j})dy+2(\psi_1+\phi \mathbf{i})dz.
\end{equation*}
By direct computation, we obtain
\begin{equation*}
dU=D\overline{w}\overline{d\sigma}.
\end{equation*}
If the stream lines lie on the surface $\partial U$, then $d\psi_1=d\psi_2=d\psi_3=0$, from which it follows that $dU=d\overline{U}$. Restricting to the surface, we thus have
\begin{equation*}
D\overline{w}\overline{d\sigma} =-d\sigma w\overline{D}.
\end{equation*}
Under these assumptions, the force integral can be expressed using integrals of the form in Proposition~\ref{ZIL}. Hence, we have:
\begin{theorem}[Monogenic Blasius–Chaplygin]\label{VOikaaava}
\begin{equation*}
F= -\frac{\rho}{8}\int_{\partial U}\Big( \text{Sc}\big((w\overline{D}) d\sigma (w\overline{D})\big)+\text{Sc}\big(  (w\overline{D}) d\sigma (w\overline{D}) \mathbf{i}\big) \mathbf{i} +  \text{Sc}\big((w\overline{D}) d\sigma (w\overline{D})\mathbf{j}\big) \mathbf{j}  \Big).
\end{equation*}
\end{theorem}
This result has an important consequence, namely:
\begin{cor}
If $w\colon U\to \mathbb{H}$ is monogenic, then $F=0$.
\end{cor}
\begin{proof}
The result follows from Proposition~\ref{ZIL}, by noting that
$w\overline{D}D=\Delta w=0$ and $D(wD)=(Dw)D=0$.
\end{proof}
It is necessary to make the following remark:
\begin{rem}
The idea in the proof of the previous theorem is to restrict the Cauchy–Riemann operator to the surface, using suitably chosen coordinates. This idea can be developed considerably further, and we refer to \cite{OS} for further details.
\end{rem}\par

\subsection*{2D reduction of the Blasius–Chaplygin force formula}

In this section, we illustrate how the classical two-dimensional Blasius–Chaplygin force formula can be obtained from the monogenic formula introduced in the previous section. Let us assume that the flow is invariant in the $z$-direction, that is, the velocity potential is $\phi = \phi(x,y)$. Then the monogenic flow potential is also invariant in the $z$-direction, and the Cauchy–Riemann operator reduces to
\begin{equation*}
D_{2d} = \frac{\partial}{\partial x} + \mathbf{i} \frac{\partial}{\partial y}.
\end{equation*}
The monogenicity condition for the flow potential then reduces to
\begin{equation*}
D_{2d}w = D_{2d}(\phi + \psi_1 \mathbf{i}) + D_{2d}(\psi_2 \mathbf{j} + \psi_3 \mathbf{k}) = 0,
\end{equation*}
which holds if equations
\begin{equation*}
D_{2d}(\phi + \psi_1 \mathbf{i}) = \partial_x \phi - \partial_y \psi_1 + (\partial_y \phi + \partial_x \psi_1)\mathbf{i} = 0
\end{equation*}
and
\begin{equation*}
D_{2d}(\psi_2 \mathbf{j} + \psi_3 \mathbf{k}) = (\partial_x \psi_2 - \partial_y \psi_3)\mathbf{j} + (\partial_x \psi_3 + \partial_y \psi_2)\mathbf{k} = 0
\end{equation*}
are satisfied.\par
As it can clearly seen from the equations presented above, the functions $\psi_2$ and $\psi_3$ are independent of the velocity potential $\phi$. This is because, as mentioned previously, a vector-valued monogenic function can be added to the stream function, that is, $\Psi \mapsto \Psi + \mathbf{H}$, and this monogenic function is naturally independent of the velocity potential. Consequently, in the above, $\mathbf{H} = \psi_2 \mathbf{j} + \psi_3 \mathbf{k}$, and we may consider the potential
\begin{equation*}
w = \phi + \psi_1 \mathbf{i},
\end{equation*}
which is precisely the two-dimensional analytic flow potential that satisfies conditions~(\ref{nop}).\par
Assume now that $C$ is a piecewise smooth curve in the plane with positive orientation (counterclockwise), that is, if $n = (n_1, n_2)$ is the outer unit normal and
\begin{equation}\label{TASONORMAALI}
(n_1 + \mathbf{i}\, n_2)ds = -\mathbf{i} \, d\mathsf{z},  
\end{equation}
where $d\mathsf{z} = dx + \mathbf{i}\,dy$. The curve $C$ defines in the space $\mathbb{R}^3$ a cylinder $C \times [-h, h]$, where $h > 0$. Let $U_h$ be a regular domain consisting of the aforementioned cylinder together with its top and bottom surfaces, that is,
\begin{equation*}
U_h = C \times [-h, h] \cup U_{\text{top}} \cup U_{\text{bottom}}.
\end{equation*}
Next, we apply formula~(\ref{BlasiusVoima}), which gives
\begin{equation*}
F=\frac{\rho}{8}\int_{C \times [-h, h]} |D_{2d}\overline{w}|^2 \, d\sigma+\frac{\rho}{8}\int_{ U_{\text{top}}} |D_{2d}\overline{w}|^2 \, d\sigma+\frac{\rho}{8}\int_{U_{\text{bottom}}} |D_{2d}\overline{w}|^2 \, d\sigma.
\end{equation*}
Since the integrand $|D_{2d}\overline{w}|^2$ does not vary in the $z$-direction, and by using~(\ref{normaali}), we have $d\sigma = \mathbf{j}dS$ on the top surface and $d\sigma = -\mathbf{j}dS$ on the bottom surface. Therefore, the last two integrals cancel each other out, since the integration goes over the same region. Thus,
\begin{equation*}
F=\frac{\rho}{8}\int_{C \times [-h, h]} |D_{2d}\BAR{w}|^2 \, d\sigma.
\end{equation*}\par
Let now the curve $C$ be parametrized by arc length as $(x(s), y(s))$ for $s \in [0, L]$, and then the cylindrical surface is parametrized as
\begin{equation*}
r(s,t) = (x(s), y(s), t),
\end{equation*}
with $t \in [-h, h]$. Then we have
\begin{equation*}
dx = x'(s)\, ds, \qquad dy = y'(s)\, ds, \qquad dz = dt,
\end{equation*}
and the surface element~(\ref{pintaalkio}) becomes
\begin{equation*}
d\sigma|_{C \times [-h, h]} = (y'(s) - \mathbf{i} x'(s))\, ds \wedge dt.
\end{equation*}
Since the tangent vector of the curve is $(x'(s), y'(s))$, the outward normal vector is obtained by rotating this vector $\frac{\pi}{2}$ clockwise, that is, $(y'(s), -x'(s))$. In terms of complex numbers, rotation corresponds to multiplication by the number $-\mathbf{i}$, hence
\begin{equation*}
n(s) = -\mathbf{i} (x'(s) + \mathbf{i} y'(s)) = y'(s) - \mathbf{i} x'(s).
\end{equation*}
This is a unit normal vector because the parametrization is by arc length. Using~(\ref{TASONORMAALI}), we can rewrite the integral as
\begin{equation*}
F = -\frac{\mathbf{i}\rho}{8} \int_{-h}^h \oint_C |D_{2d}\BAR{w}|^2\, d\mathsf{z}\, dt.
\end{equation*}\par
To obtain a two-dimensional formula for the force density, we define
\begin{equation*}
F_{2d} := \frac{F}{2h} = -\frac{\mathbf{i}\rho}{16h} \int_{-h}^h dt \oint_C |D_{2d}\overline{w}|^2 \, d\mathsf{z} = -\frac{\mathbf{i}\rho}{8} \oint_C |D_{2d}\overline{w}|^2 \, d\mathsf{z}.
\end{equation*}
Let now $f(\mathsf{z}) = \phi + \mathbf{i} \psi_1$, so that $f'(\mathsf{z}) = v_1 - \mathbf{i} v_2$. Then
\begin{equation*}
D_{2d}\overline{w}=2 \overline{f'(\mathsf{z})},
\end{equation*}
and thus
\begin{equation*}
F_{2d} = -\frac{\mathbf{i}\rho}{2} \oint_C |f'(\mathsf{z})|^2 \, d\mathsf{z},
\end{equation*}
where $F_{2d} = F_x + \mathbf{i} F_y$. From this point onward, the derivation of the formula is a standard task in fluid mechanics textbooks. That is, by taking the complex conjugate we obtain
\begin{equation*}
F_x - \mathbf{i}F_y=\frac{\mathbf{i}\rho}{2} \oint_C f'(\mathsf{z})\overline{f'(\mathsf{z})} \, d\overline{\mathsf{z}}=\frac{\mathbf{i}\rho}{2} \oint_C f'(\mathsf{z})\, d\overline{f},
\end{equation*}
where $df=d\phi+\mathbf{i}d\psi_1$. Since $d\psi_1=0$ along the streamlines, we have $d\overline{f}=df=f'(\mathsf{z})  d\mathsf{z}$, hence
\begin{equation*}
F_x - \mathbf{i} F_y=\frac{\mathbf{i}\rho}{2} \oint_C (f'(\mathsf{z}))^2\, d\mathsf{z},
\end{equation*}
which is the holomorphic form of the Blasius–Chaplygin formula.\par

\section{Quaternionic Blasius–Chaplygin formula for the moment}

In this section, we introduce a quaternionic Blasius–Chaplygin formula for the moment. Assume that $U$ is a solid body and that the pressure field $p(x,y,z)$ is known on its surface. Suppose that $\mathsf{x}_0 \in \mathbb{R}^3$ is an arbitrary point. We work with reduced quaternions and define the moment by
\begin{equation*}
M := r\times F=\begin{vmatrix} 1 & \mathbf{i} & \mathbf{j} \\ r_1 & r_2 &  r_3\\ F_1 & F_2 & F_3 \end{vmatrix} 
=r_2F_3-r_3F_2+(r_3F_1-r_1F_3)\mathbf{i}+(r_1F_2-r_2F_1)\mathbf{j},
\end{equation*}
where $r = r_1 + r_2\mathbf{i} + r_3\mathbf{j}$ and $F = F_1 + F_2\mathbf{i} + F_3\mathbf{j}$.\par
Let $\mathsf{x} \in \partial U$ and set $r(\mathsf{x}) = \mathsf{x} - \mathsf{x}_0$, then the moment induced by the pressure field with respect to the point $\mathsf{x}_0$ is given by
\begin{equation*}
M_{\mathsf{x}_0} = -\int_{\partial U} p(\mathsf{x})(r(\mathsf{x}) \times n(\mathsf{x})) dS(\mathsf{x}),
\end{equation*}
where the cross product is interpreted in the sense of reduced quaternions, as in the definition of the moment presented above.\par 
From the definition of the moment, we immediately obtain the following proposition:
\begin{prop}\label{ANKO}
If the pressure $p$ is constant on the entire surface $\partial U$, then $M_{\mathsf{x}_0} = 0$.
\end{prop}
\begin{proof}
The moment can alternatively be computed as
\begin{equation*}
M_{\mathsf{x}_0} = (\mathsf{x}_c - \mathsf{x}_0) \times F,
\end{equation*}
where $\mathsf{x}_c$ is the center of mass of the body, and $F$ is the total force as computed in the previous result. When the pressure is constant, then $F = 0$.
\end{proof}\par
Substituting Bernoulli’s formula, we obtain
\begin{equation*}
M_{\mathsf{x}_0} = \frac{\rho}{2} \int_{\partial U} |v(\mathsf{x})|^2 (r(\mathsf{x}) \times n(\mathsf{x})) dS(\mathsf{x}),
\end{equation*}
leading to the \textit{Blasius–Chaplygin formula for the moment}
\begin{equation*}
M_{\mathsf{x}_0} = \frac{\rho}{8} \int_{\partial U} |D\overline{w}(\mathsf{x})|^2 r(\mathsf{x}) \times d\sigma(\mathsf{x}).
\end{equation*}
The force formula was presented in Theorem \ref{VOikaaava} using monogenic functions. For the moment, the situation is more complicated, and the formula is left as an open problem in this article. However, by using the same proof technique as in Proposition \ref{ANKO}, it can be shown that:

\begin{prop}
If $w:U\to \mathbb{H}$ is monogenic, then $M_{\mathsf{x}_0} = 0$.
\end{prop}

\subsection*{2D reduction of the Blasius–Chaplygin moment formula}

Finally, in this section, we show how to obtain the classical two-dimensional Blasius–Chaplygin moment formula from the monogenic formula introduced in the previous section. Assume $\mathsf{x}_0 = 0$, in which case
\begin{equation*}
M = \frac{\rho}{8} \int_{\partial U} |D\overline{w}(\mathsf{x})|^2 (\mathsf{x} \times n(\mathsf{x})) dS(\mathsf{x}).
\end{equation*}
Again, suppose the flow is invariant in the $z$-direction, and write the integral over the previously defined surface $U_h$, which leads to
\begin{equation*}
M = \frac{\rho}{8} \int_{-h}^h \int_0^L |D_{2d}\overline{w}|^2 (\mathsf{x}(s,t) \times n(s)) \, ds \, dt, 
\end{equation*}
where
\begin{equation*}
\mathsf{x}(s,t) = x(s) + y(s)\mathbf{i} + t \mathbf{j}
\end{equation*}
and
\begin{equation*}
n(s) = y'(s) - \mathbf{i}x'(s),
\end{equation*}
implying that
\begin{equation*}
\mathsf{x}(s,t) \times n(s) = (x'(s) + y'(s)\mathbf{i})t - (x(s)x'(s) + y(s)y'(s))\mathbf{j}.
\end{equation*}
Substituting this expression into the definition of $M$ and performing the outer integration yields
\begin{equation*}
M = -\mathbf{j} \frac{\rho h}{4} \int_0^L |D_{2d}\overline{w}|^2 (x(s)x'(s) + y(s)y'(s)) \, ds.
\end{equation*}
Let $D_{2d}\overline{w} = 2 \overline{f'(\mathsf{z})}$, so we obtain
\begin{equation*}
M = -\mathbf{j} \rho h \int_0^L |f'(x(s), y(s))|^2 (x(s)x'(s) + y(s)y'(s)) \, ds.
\end{equation*}
Note that $dx = x'(s) ds$ and $dy = y'(s) ds$, meaning that $\text{Re}(\mathsf{z} d\overline{\mathsf{z}}) = (x(s)x'(s) + y(s)y'(s)) ds$. Consequently, the moment density becomes
\begin{equation*}
M_{2d} = \frac{M}{2h} = -\frac{\rho}{2} \text{Re} \left( \oint_C \mathsf{z} |f'(\mathsf{z})|^2 \, d\overline{\mathsf{z}} \right) \mathbf{j}.
\end{equation*}
Proceeding analogously as with the force formula, we finally obtain
\begin{equation*}
M_{2d} = -\frac{\rho}{2} \text{Re} \left( \oint_C \mathsf{z} (f'(\mathsf{z}))^2 \, d\mathsf{z} \right) \mathbf{j},
\end{equation*}
which is the two-dimensional Blasius–Chaplygin formula for the moment around the $\mathbf{j}$-axis.\\

Finally, it should be noted that, unlike in two dimensions, the three-dimensional Blasius–Chaplygin formulas yield actual forces and moments instead of force and moment distributions. However, the assumption of two-dimensionality is always an idealization of the physical problem.

\section*{Conclusion}
The main result in this article is the proof of a Blasius–Chaplygin type force formula in three dimensions when the flow potential is quaternion-valued. This result is presented in Theorem \ref{VOikaaava}. It demonstrates that quaternion analysis can be used in the study of 3D ideal flow in a manner analogous to complex analysis, and in fact, the complex 2D idealization is always obtained as a special case. The results presented in this article open up possibilities for the investigation of three-dimensional fluid dynamics problems, see for example \cite{LS}. \\

From the perspective of the mathematical tools employed, quaternion algebra is used instead of complex numbers, with the most significant difference being the non-commutativity of multiplication. Moreover, instead of contour integrals, the integrals are surface integrals in the space $\mathbb{R}^3$, whose effective handling requires the use of differential forms. The aforementioned mathematics demands a somewhat different mindset from the practitioner, although all concepts and tools from complex analysis have their counterparts in quaternion analysis.

\end{document}